\documentclass[12pt]{article}
\usepackage{amsmath}
\usepackage{amssymb}
\usepackage{amsfonts}
\usepackage{amscd}
\usepackage{enumerate}
\usepackage{fancyheadings}
\usepackage{enumerate}
\usepackage{color}
\usepackage{graphicx}
\newtheorem{Lem}{Lemma}

\newtheorem{Th}{Theorem}

\newtheorem{Def}{Definition}

\DeclareMathOperator{\supp}{supp}

\DeclareMathOperator{\esupp}{ess-supp}

\begin{document}
\begin{center}
\bf{\Large Semi-Classical Wavefront Set and Fourier Integral Operators}
\end{center}

\noindent
{\bf Ivana Alexandrova}

\noindent
Department of Mathematics, University of Toronto, Toronto, Ontario, Canada M5S 3G3, alexandr@math.toronto.edu

\noindent
July 26, 2004

\begin{abstract}
Here we define and prove some properties of the semi-classical wavefront set.
We also define and study semi-classical Fourier integral operators and prove a generalization of Egorov's Theorem to manifolds of different dimensions. 
\end{abstract}


%\newpage
\section{Introduction}
In this article we define and establish some of the properties of the semi-classical wavefront set and semi-classical Fourier integral operators.
The paper is organized as follows.
In Section \ref{sprelims} we review some of the theory of semi-classical pseudodifferential operators, which we will use here.
In Section \ref{shWF} we gather the existing definitions of semi-classical wavefront set and show that they are equivalent.
We furthe explore the properties of the semi-classical wavefront set in Section \ref{sprops}.
In Section \ref{shFIO} we define and prove a characterization of global semi-classical Fourier integral operators as well as a generalization of Egorov's Theorem to manifolds of unequal dimensions.

\section{Preliminaries}\label{sprelims}
In this section we recall some of the elements of semi-classical analysis 
which we 
will use here.
First we define two classes of symbols
\begin{equation*}
S_{2n}^{m}\left(1\right)= \left\{ a\in
C^{\infty}\left(\mathbb{R}^{2n}\times(0, h_0]\right): \forall
\alpha, \beta\in\mathbb{N}^{n}, \sup_{(x, \xi, 
h)\in\mathbb{R}^{2n}\times (0, 
h_{0}]}h^{m}\left|\partial^{\alpha}_{x}\partial^{\beta}_{\xi}a\left(x, 
\xi;
h\right)\right|\leq
C_{\alpha, \beta}\right\}
\end{equation*}
and
\begin{equation*}
S^{m, k}\left(T^{*}\mathbb{R}^{n}\right)=\left\{a\in
C^{\infty}\left(T^{*}\mathbb{R}^{n}\times(0, h_0]\right): \forall \alpha, 
\beta\in\mathbb{N}^{n}, \left|\partial^{\alpha}_{x}\partial^{\beta}_{\xi} 
a\left(x, \xi;
h\right)\right|\leq
C_{\alpha,
\beta}h^{-m}\left\langle\xi\right\rangle^{k-|\beta|}\right\},
\end{equation*}
where $h_0\in(0,1]$ and $m, k\in\mathbb{R}.$
For $a\in S_{2n}\left(1\right)$ or $a\in S^{m, k}\left(T^{*}\mathbb{R}^{n}\right)$ we define the
corresponding semi-classical pseudodifferential operator of class $\Psi_{h, t}^{m}(1, \mathbb{R}^{n})$ or $\Psi_{h, t}^{m, k}(\mathbb{R}^{n}),$ respectively, by setting
\begin{equation*}
Op_{h}\left(a\right)u\left(x\right)=\frac{1}{\left(2\pi
h\right)^{n}}\int\int e^{\frac{i\left\langle x-y,
\xi\right\rangle}{h}}a\left(x, \xi; h\right)u\left(y\right) dy d\xi, \;u\in
\mathcal{S}\left(\mathbb{R}^{n}\right)
\end{equation*}
for $t\in [0, 1]$ and extending the definition to $\mathcal{S}'\left(\mathbb{R}^{n}\right)$ by
duality
(see \cite{DS}).
Below we shall work only with symbols which admit asymptotic expansions in $h$ and with
pseudodifferential operators which are quantizations of such symbols.
For $A\in\Psi_{h, t}^{k}(1, \mathbb{R}^{n})$ or $A\in\Psi_{h, t}^{m, k}(\mathbb{R}^{n}),$ we shall use $\sigma_{0}(A)$ 
and $\sigma(A)$ to denote its principal symbol and its complete symbol, respectively.

For $a\in S^{m, k}_{n}\left(T^{*}\mathbb{R}^{n}\right)$ we
define:
\begin{equation*}
\begin{aligned}
& \esupp_{h} a\\
&\quad =\Big\{\left(x, \xi\right)\in T^{*}X|\: \exists\:
\epsilon>0\;
\partial_{x}^{\alpha}\partial_{\xi}^{\beta}a\left(x',
\xi'\right)=\mathcal{O}_{C(B((x, \xi),
\epsilon))}\left(h^{\infty}\right),\; \forall\alpha,
\beta\in\mathbb{N}^{n}\Big\}^{c}\\
& \quad\cup\bigg(\bigg\{\left(x, \xi\right)\in
T^{*}X\backslash\left\{0\right\}|\:
\exists\:\epsilon>0\:
\partial_{x}^{\alpha}\partial_{\xi}^{\beta}a\left(x',
\xi'\right)=\mathcal{O}\left(h^{\infty}\left\langle\xi\right\rangle^{-\infty}\right),\\
&\quad\quad \text{uniformly in } (x', \xi') \text{ such that }
\|x-x'\|+\frac{1}{|\xi'|}+\left|\frac{\xi}{|\xi|}-\frac{\xi'}{|\xi'|}\right|<\epsilon\bigg\}/
\mathbb{R}_{+}\bigg)^{c}\\
&\quad \subset T^{*}X\sqcup S^{*}X,
\end{aligned}
\end{equation*}
where we define
$S^{*}X=\left(T^{*}X\backslash\left\{0\right\}\right)/\mathbb{R}_{+}.$
For $A\in\Psi^{m, k}_{h}\left(\mathbb{R}^{n}\right),$ we then define
\begin{equation*}
WF_{h}\left(A\right)=\esupp_{h} a, A=Op_{h}\left(a\right).
\end{equation*}

We also define the class of semi-classical distributions 
$\mathcal{D}_{h}'(\mathbb{R}^{n})$ with which we will work here
\begin{equation*}
\begin{aligned}
\mathcal{D}'_{h}(\mathbb{R}^{n}) = & \big\{u\in C^{\infty}_{h}\left((0, 
1];
\mathcal{D}'\left(\mathbb{R}^{n}\right)\right): \forall\chi\in 
C_{c}^{\infty}\left(\mathbb{R}^{n}\right) \exists\: N\in\mathbb{N}\text{ 
and 
} C_{N}>0:\\
& \quad |\mathcal{F}_{h}\left(\chi u\right)\left(\xi, h\right)|\leq
C_{N}h^{-N}\langle\xi\rangle^{N}\big\}
\end{aligned}
\end{equation*}
where
\begin{equation*}
\mathcal{F}_{h}\left(u\right)\left(\xi,
h\right)=\int_{\mathbb{R}^{n}}e^{-\frac{i}{h}\left\langle x,
\xi\right\rangle}u\left(x, h\right)dx
\end{equation*}
with the obvious extension of this definition to 
$\mathcal{E}_{h}'(\mathbb{R}^{n}).$
We shall work with the $L^{2}-$ based semi-classical Sobolev spaces 
$H^{s}(\mathbb{R}^{n}),$ $s\in\mathbb{R},$ which consist of the distributions 
$u\in\mathcal{D}_{h}'(\mathbb{R}^{n})$ such 
that 
$\|u\|_{H^{s}(\mathbb{R}^{n})}^{2}=\frac{1}{(2\pi 
h)^{n}}\int_{\mathbb{R}^{n}}(1+\|\xi\|^{2})^{s}\left|\mathcal{F}_{h}(u)(\xi, 
h)\right|^{2}d\xi<\infty.$

We shall say that $u=v$ {\it microlocally} near an open set
$U\subset T^{*}\mathbb{R}^{n}$, if
$P(u-v)=\mathcal{O}\left(h^{\infty}\right)$ in
$C_{c}^{\infty}\left(\mathbb{R}^{n}\right)$ for
every $P\in \Psi^{0}_{h}\left(1, \mathbb{R}^{n}\right)$ such that
\begin{equation}\label{P}
WF_{h}\left(P\right)\subset \tilde{U}, \bar{U}\Subset \tilde{U}\Subset
T^{*}\mathbb{R}^{n}, \tilde{U} \text{ open}.
\end{equation}
We shall also say that $u$ satisfies a property $\mathcal{P}$  {\it
microlocally} near an open set $U\subset T^{*}{\mathbb{R}^{n}}$ if there
exists $v\in\mathcal{D}_{h}'\left(\mathbb{R}^{n}\right)$ such that $u=v$
microlocally
near $U$ and $v$ satisfies property $\mathcal{P}$.

For open sets $U, V\subset T^{*}\mathbb{R}^{n},$ the operators $T,
T'\in\Psi^{m}_{h}\left(\mathbb{R}^{n}\right)$ are said to be {\it
microlocally
equivalent} near $V\times U$ if for any $A, 
B\in\Psi_{h}^{0}\left(\mathbb{R}^{n}\right)$
such that
\begin{equation*}
WF_{h}\left(A\right)\subset\tilde{V},
WF_{h}\left(B\right)\subset\tilde{U},
\bar{V}\Subset\tilde{V}\Subset T^{*}\mathbb{R}^{n},
\bar{U}\Subset\tilde{U}\Subset T^{*}\mathbb{R}^{n}, \tilde{U}, \tilde{V}
\text{ open }
\end{equation*}
\begin{equation*}
A\left(T-T'\right)B=\mathcal{O}\left(h^{\infty}\right)\colon\mathcal{D}_{h}'\left(\mathbb{R}^{n}\right)\rightarrow
C^{\infty}\left(\mathbb{R}^{n}\right).
\end{equation*}
We shall also use the notation $T\equiv T'.$

\section{Semi-Classical Wavefront Set}\label{shWF}

In this section we discuss the different notions of semi-classical 
wavefront set used in the 
literature and show that they are equivalent.
We further establish some of their properties.

We further let $\hat{T}^{*}\mathbb{R}^{n}=T^{*}\mathbb{R}^{n}\sqcup 
S^{*}\mathbb{R}^{n},$ where we set 
$S^{*}\mathbb{R}^{n}=\left(T^{*}\mathbb{R}^{n}\backslash 0 
\right)/\mathbb{R}_{+}$ with the $\mathbb{R}_{+}$ action given by 
mutiplication on the fibers: $(x, \xi)\mapsto (x, t\xi).$
As in \cite{Michel}, the points in $T^{*}\mathbb{R}^{n}$ will be called 
finite 
and the points in 
$S^{*}\mathbb{R}^{n}$ will be called infinite.
We make the following 
definition as in \cite{Michel}
\begin{Def}\label{defM}
Let $u\in\mathcal{D}'_{h}\left(\mathbb{R}^{n}\right)$ and let 
$\left(x_{0}, 
\xi_{0}\right)\in\hat{T}^{*}\left(\mathbb{R}^{n}\right).$
We shall say that $\left(x_{0}, \xi_{0}\right)$ does not belong to 
$WF_{h}\left(u\right)$ if:
\begin{itemize}
\item If $\left(x_{0}, \xi_{0}\right)$ is finite: there exist $\chi\in 
C_{c}^{\infty}\left(\mathbb{R}^{n}\right)$ with $\chi\left(x_{0}\right)\ne 
0$ and an open neighborhood 
$U$ of $\xi_{0}$, such that $\forall N\in\mathbb{N},$ $\forall\xi\in U,$ 
$|\mathcal{F}\left(\chi u\right)\left(\xi, h\right)|\leq C_{N}h^{N}.$
We shall denote the complement of the set of all such points by $WF_{h}^{f}(u).$
\item If $\left(x_{0}, \xi_{0}\right)$ is infinite: there exist $\chi\in
C_{c}^{\infty}\left(\mathbb{R}^{n}\right)$ with $\chi\left(x_{0}\right)\ne 
0$ and a 
conic neighborhood $U$ of $\xi_{0}$, such that $\forall N\in\mathbb{N},$ 
$\forall\xi\in 
U\cap\left\{|\xi|\geq\frac{1}{C}\right\},$
\[|\mathcal{F}\left(\chi u\right)\left(\xi, h\right)|\leq 
C_{N}h^{N}\left\langle\xi\right\rangle^{-N}.\]
We shall denote the complement of the set of all such points by $WF_{h}^{i}(u).$
\end{itemize}
\end{Def}

The definition of semi-classical wavefront set given in \cite{SZQ} is as 
follows
\begin{Def}\label{defSM}
\[
WF_{h}\left(u\right)=\left\{\left(x, \xi\right): \exists A\in\Psi_{h}^{0, 
0}\left(\mathbb{R}^{n}\right) 
\sigma_{o}\left(A\right)\left(x, \xi\right)\ne 0,\\
 A u\in h^{\infty}C^{\infty}\left((0, 1]_{h}; 
C^{\infty}\left(\mathbb{R}^{n}\right)\right)\right\}^{c}.\]
\end{Def}

\begin{Lem}
$u\in\mathcal{D}'_{h}\left(\mathbb{R}^{n}\right)$ if and only if for every 
$\chi\in 
C_{c}^{\infty}\left(\mathbb{R}^{n}\right)$ there exist
$m, k_m, C_{m} \in\mathbb{R}$ such that for every  
$\|\chi u\|_{H^{m}\left(\mathbb{R}^{n}\right)}\leq C_{m}h^{-k_{m}}.$
\end{Lem}

\medskip
\noindent
{\it Proof:}\;
The first implication is clear.
For the second implication, let $s\in\mathbb{R}$ be such that 
$m+s>\frac{n}{2}$ and 
consider
\begin{equation*}
\begin{aligned}
& C_{m}^{2}h^{-2k_{m}}\geq \frac{1}{\left(2\pi h\right)^{n}}\int 
\left(1+|\xi|^{2}\right)^{m}|\mathcal{F}_{h}\left(\chi 
u\right)\left(\xi, h\right)|^{2}d\xi\\
& = \frac{1}{\left(2\pi h\right)^{n}}\int 
\left(1+|\xi|^{2}\right)^{m+s}\frac{|\mathcal{F}_{h}\left(\chi
u\right)\left(\xi, 
h\right)|^{2}}{\left(1+|\xi|^{2}\right)^{s}}d\xi
=\|\tilde{u}\|^{2}_{H^{m+s}\left(\mathbb{R}^{n}\right)},
\end{aligned}
\end{equation*}
where 
$\tilde{u}=\mathcal{F}_{h}^{-1}\left(\frac{|\mathcal{F}_{h}\left(\chi 
u\right)\left(\cdot, h\right)}
{\left(1+|\cdot|^{2}\right)^{\frac{s}{2}}}\right).$
Then we have that
\begin{equation*}
C_{m}h^{-k_{m}}\geq \|\tilde{u}\|_{H^{m+s}\left(\mathbb{R}^{n}\right)}\geq
C\|\tilde{u}\|_{L^{\infty}\left(\mathbb{R}^{n}\right)}\geq
C\frac{1}{|\supp 
\chi|}\left\|\hat{\tilde{u}}\right\|_{L^{\infty}(\mathbb{R}^{n})}
\end{equation*}
and therefore
\begin{equation*}
|\mathcal{F}_{h}\left(\chi u\right)\left(\xi, h\right)|\leq 
C h^{-k_m} (1+|\xi|^{2})^{\frac{s}{2}}.
\end{equation*}$\hfill\Box$

\begin{Lem}
Definitions (\ref{defM}) and (\ref{defSM}) are 
equivalent.
\end{Lem}

\medskip
\noindent
{\it Proof:}\;
Let $\left(x_{0}, \xi_{0}\right)\in 
T^{*}\mathbb{R}^{n}\backslash 
WF_{h}^{f}\left(u\right).$
Let $\varphi\in C_{c}^{\infty}\left(\mathbb{R}^{n}\right)$ satisfy 
$\varphi\left(x_{0}\right)\ne 0$ and let
$\chi\in C_{c}^{\infty}\left(\mathbb{R}^{n}\right)$ have support in a 
bounded open neighborhood $V$ of $\xi_{0}$ such that
$\mathcal{F}_{h}\left(\varphi 
u\right)\left(\xi\right)=\mathcal{O}\left(h^{\infty}\right)$ uniformly for 
$\xi \in V.$
Consider
\begin{equation*}
Au\left(x\right)=\frac{1}{\left(2\pi h\right)^{n}}\int\int 
e^{\frac{i}{h}\left\langle x-y,
\xi\right\rangle}\varphi\left(x\right)\varphi\left(y\right)\chi\left(\xi\right)u\left(y\right) 
d y d \xi.
\end{equation*}
We clearly have that $A\in \Psi_{h}^{0, 0}\left(\mathbb{R}^{n}\right)$ 
with
$\sigma_{0}^{l}\left(A\right)\left(x_{0}, \xi_{0}\right)\ne 0$ and 
$Au=\mathcal{O}\left(h^{\infty}\right)$ in
$C^{\infty}\left(\mathbb{R}^{n}\right),$ where $\sigma^{l}_{0}\left(A\right)$ 
is the principal symbol of the
left-quantization of $A.$

Let, now, $\left(x_{0}, \xi_{0}\right)\in T^{*}\mathbb{R}^{n}$ be such 
that 
there exists
$A\in\Psi_{h}^{0, 0}\left(\mathbb{R}^{n}\right)$ elliptic at $\left(x_{0}, 
\xi_{0}\right)$ such that
$Au=\mathcal{O}\left(h^{\infty}\right)$ in 
$C^{\infty}\left(\mathbb{R}^{n}\right).$
Let $\varphi, \chi\in C_{c}^{\infty}\left(\mathbb{R}^{n}\right)$ be
such that $\varphi\left(x_{0}\right)\ne
0,$ $\chi\left(\xi_{0}\right)\ne 0,$ and 
$\chi\left(hD\right)\varphi=BA+R,$ where $B\in\Psi_{h}^{0, 
0}\left(\mathbb{R}^{n}\right),$ $R\in \Psi_{h}^{-\infty, 
-\infty}\left(\mathbb{R}^{n}\right).$
Then
\begin{equation*}
\begin{aligned}
\chi\left(hD\right)\varphi u\left(x\right) & =\frac{1}{\left(2\pi 
h\right)^{n}}\int\int
e^{\frac{i}{h}\left\langle x-y,
\xi\right\rangle}\chi\left(\xi\right)\varphi\left(y\right)u\left(y\right)dy d\xi\\
& = \frac{1}{\left(2\pi h\right)^{n}}\int e^{\frac{i}{h} x\cdot 
\xi}\chi\left(\xi\right)\widehat{\varphi u}\left(\xi/h\right) d\xi
=\mathcal{O}\left(h^{\infty}\right) \text{ in } 
C^{\infty}\left(\mathbb{R}^{n}\right).
\end{aligned}
\end{equation*}
Therefore, $\chi\left(\xi\right)\widehat{\varphi 
u}\left(\xi/h\right)=\mathcal{O}\left(h^{\infty}\right)$
uniformly in $\xi$ and therefore $\widehat{\varphi
u}\left(\xi/h\right)=\mathcal{O}\left(h^{\infty}\right)$ uniformly in 
$\xi$ in a bounded open
set containing $\xi_{0},$ which implies that $\left(x_{0}, 
\xi_{0}\right)\notin
WF_{h}^{f}\left(u\right).$

The case of an infinite point is handled similarly.
See also \cite{GS} for the proof in the classical setting, which applies directly to the 
infinite semi-classical wavefront set here.  $\hfill\Box$

\subsection{Properties of the Semi-classical Wavefront Set}\label{sprops}

In this section, we prove the following properties of the semi-classical wavefront set
\begin{Lem}\label{wfhprop}
Let $u\in\mathcal{D}_{h}'(\mathbb{R}^{d_{1}}),$ 
$v\in\mathcal{D}_{h}'(\mathbb{R}^{d_{2}}),$ 
$w\in\mathcal{D}_{h}'(\mathbb{R}^{d_{3}}),$
$V\in\mathcal{D}_{h}'(\mathbb{R}^{d_{1}+d_{2}}),$ 
$W\in\mathcal{D}_{h}'(\mathbb{R}^{d_{2}+d_{3}}).$
Then
%\begin{enumerate} 
%\noindent
%\item
\begin{enumerate}[(a)]
\item $WF_{h}^{f}\left(Au\right)\subset WF_{h}^{f}\left(A\right)\cap 
WF_{h}^{f}\left(u\right)$ and $WF_{h}^{i}\left(Au\right)\subset 
WF_{h}^{i}\left(A\right)\cap WF_{h}^{i}\left(u\right)$, for $A\in
\Psi_{h}^{m, k}\left(\mathbb{R}^{d_{1}}\right).$
\item $u\otimes v\in\mathcal{D}_{h}'(\mathbb{R}^{d_{1}+d_{2}}),$ 
$WF_{h}^{f}\left(u\otimes v\right)\subset 
WF_{h}^{f}\left(u\right)\times 
WF_{h}^{f}\left(v\right)$ and 
\[\begin{aligned}
WF_{h}^{i}\left(u\otimes 
v\right)\subset & \left(WF_{h}^{i}\left(u\right)\times
WF_{h}^{i}\left(v\right)\right)\cup
\left(\left(\supp u \times\{0\}\right)\times 
WF_{h}^{i}\left(v\right)\right)\\
& \quad \cup\left(WF_{h}^{i}\left(u\right)\times\left(\supp 
v\times\{0\}\right)\right).
\end{aligned}\]
\item if $V$ is proper,
$WF_{h}^{i}(v)\cap\left(WF_{h}^{i}\right)'_{\mathbb{R}^{d_{2}}}(V)=\emptyset,$
and $WF_{h}^{f}(v)$
is compact, then $Vv\in\mathcal{D}_{h}'(\mathbb{R}^{d_{1}})$ and
\[WF_{h}^{i}\left(V
v\right)\subset
\left(WF_{h}^{i}\right)'\left(V\right)\left(WF_{h}^{i}\left(v\right)\right)
\cup\left(WF_{h}^{i}\right)'_{\mathbb{R}^{d_{1}}}(V)\]
and 
\[WF_{h}^{f}\left(V v\right)\subset
\left(WF_{h}^{f}\right)'\left(V\right)\left(WF_{h}^{f}\left(v\right)\right),\] 
where
\[
\left(WF_{h}^{i}\right)'_{\mathbb{R}^{d_{2}}}(V)=\Big\{(y, 
\eta)\in
T^{*}\mathbb{R}^{d_{2}} \backslash
\{0\}: \exists x\in\mathbb{R}^{d_{1}}, (x, 0; y, 
\eta)\in\left(WF_{h}^{i}\right)'(V)\Big\}\] and $Vv$ is defined as in Theorem 
7.8, \cite{GS}.
The same conclusion holds if $V$ is not necessarily proper but
$v\in\mathcal{E}_{h}'(\mathbb{R}^{d_{2}})$ and all the other assumptions 
are satisfied.
\item if at least one of $V$ and $W$ are properly supported, 
\[\left(WF_{h}^{i}\right)'_{\mathbb{R}^{d_{2}}}(V)\cap
\left(WF_{h}^{i}\right)'_{\mathbb{R}^{d_{2}}}(W)=\emptyset,\]
\[\left\{p\in\mathbb{R}^{d_{2}}: \exists\: 
(q, r)\in\mathbb{R}^{d_{1}}\times\mathbb{R}^{d_{3}}, (q, p)\in 
WF_{h}^{f}(V), (p, r)\in 
WF_{h}^{f}(W)\right\}\] is compact, then $V\circ 
W\in\mathcal{D}_{h}'(\mathbb{R}^{d_{1}+d_{3}}),$ 
\[\left(WF_{h}^{i}\right)'\left(V\circ 
W\right)\subset \left(WF_{h}^{i}\right)'\left(V\right)\circ 
\left(WF_{h}^{i}\right)'\left(W\right),\] and \[\left(WF_{h}^{f}\right)'\left(V\circ 
W\right)\subset \left(WF_{h}^{f}\right)'\left(V\right)\circ 
\left(WF_{h}^{f}\right)'\left(W\right),\] where $V\circ W$ is defined as in Theorem 7.10, 
\cite{GS}.
%\end{enumerate}
\end{enumerate} 
\end{Lem}
 
{\bf Remark.}  Part (c) of this lemma is proved in \cite[Proposition 
A.I.13]{ChG} without the assumption on $WF_{h}(v).$
In our proof, however, we also show that all estimates can be made uniformly in a 
neighborhood of $WF_{h}(v).$

\medskip
\noindent
{\it Proof:}\;
In this proof we shall use $\langle\cdot, \cdot\rangle$ to 
denote 
the 
distribution pairing.
We begin by proving (a).
Let $\left(x_{0}, \xi_{0}\right)\notin WF_{h}\left(A\right)$ and 
assume that $\left(x_{0},
\xi_{0}\right)$ is a finite point.
Let $B\in\Psi_{h}^{0, 0}\left(\mathbb{R}^{d_{1}}\right)$ satisfy 
$\sigma\left(B\right)\left(x_{0},
\xi_{0}\right)\ne 0,$ $\sigma\left(B\right)\in 
C_{c}^{\infty}\left(\mathbb{R}^{d_{1}}\right),$ and
$WF_{h}\left(B\right)\cap WF_{h}\left(A\right)=\emptyset.$
Then $BA\in\Psi_{h}^{-\infty, -\infty}\left(\mathbb{R}^{d_{1}}\right)$ and 
therefore $BAu=\mathcal{O}\left(h^{\infty}\right)$ in 
$C^{\infty}\left(\mathbb{R}^{d_{1}}\right).$
If $\left(x_{0}, \xi_{0}\right)$ is an infinite point, we can again 
find
$B\in\Psi_{h}^{0, 0}\left(\mathbb{R}^{d_{1}}\right)$ with 
$WF_{h}\left(B\right)$ consisting only of
infinite points such that $WF_{h}\left(B\right)\cap 
WF_{h}\left(A\right)=\emptyset$ and we then
have that $BAu=\mathcal{O}\left(h^{\infty}\right)$ in 
$C^{\infty}\left(\mathbb{R}^{d_{1}}\right).$

Let, now, $\left(x_{0}, \xi_{0}\right)\notin WF_{h}^{f}\left(u\right).$ 
Let $c\in C_{c}^{\infty}\left(\mathbb{R}^{d_{1}}\right)$ satisfy 
$c\left(x_{0},
\xi_{0}\right)\ne 0$ and let $d\in S_{2d_{1}}\left(1\right)$ be such that 
$d\#_{h} c=1$ in a
neighborhood $W\subset \left(WF_{h}^{f}(u)\right)^{c}$ of $\left(x_{0}, 
\xi_{0}\right).$
Further, let $\chi\in C_{c}^{\infty}\left(T^{*}\mathbb{R}^{d_{1}}\right)$ have 
support in
an open set $V\Subset W$ and be equal to 1 on an open subset $U\Subset V.$
Then the operator $T=Op_{h}\left(\chi d\#_{h}c\right)$ has symbol 
$\sigma\left(T\right)\equiv
1\text{ mod } h^{\infty}$ in $S_{2d_{1}}\left(1\right)$ in $U$ and supported in $W$ and 
therefore
$Tu=\mathcal{O}\left(h^{\infty}\right)$ in 
$C^{\infty}\left(\mathbb{R}^{d_{1}}\right).$
Let $B$ be elliptic at $\left(x_{0}, \xi_{0}\right)$ with 
$WF_{h}\left(B\right)\subset U.$
Then we have that $BA\equiv BAT\text{ mod } \Psi^{-\infty, 
-\infty}_{h}\left(\mathbb{R}^{d_{1}}\right)$ and hence $BAu\equiv 
BATu=\mathcal{O}\left(h^{\infty}\right)$ in
$C^{\infty}\left(\mathbb{R}^{d_{1}}\right).$
Therefore, $\left(x_{0}, \xi_{0}\right)\notin WF_{h}\left(Au\right).$
The proof is similar in the case of an infinite point $\left(x_{0}, 
\xi_{0}\right).$

We now turn to proving (b).
It is trivial to check that $u\otimes 
v\in\mathcal{D}_{h}'(\mathbb{R}^{d_{1}+d_{2}}).$
Let $(x_0, \xi_0; y_0, 
\eta_0)\notin WF_{h}^{f}(u)\times WF_{h}^{f}(v)$ 
and let $O_1,$ $O'_{1}\subset\mathbb{R}^{d_{1}}$ and $O_2,$ 
$O'_{2}\subset\mathbb{R}^{d_{3}}$ be 
open neighbordhoods of $x_0, \xi_0, y_0,$ and $\eta_0,$ respectively, such that $O_1\times 
O'_{1}\times O_2\times O'_{2}\subset \left(WF_{h}^{f}(u)\times WF_{h}^{f}(v)\right)^{c}.$ 
Without loss of generality, we can assume that $(x_0, \xi_0)\notin WF_{h}^{f}(u).$
Then there exists $\chi_1\in C_{c}^{\infty}(\mathbb{R}^{d_{1}})$ with 
$\chi_{1}(x_0)\ne 0$ and a 
bounded open set $O'_{1}\subset\mathbb{R}^{d_{1}}$ with $\xi_0\in O'_{1}$ 
such that 
$\left|\mathcal{F}_{h}(\chi_{1}u)(\xi, h)\right|=\mathcal{O}(h^{\infty})$ uniformly for 
$\xi\in O'_{1}.$
Let, now, $\chi_2\in C_{c}^{\infty}(\mathbb{R}^{d_{2}})$ have support near 
$y_0.$ 
Then $\left|\mathcal{F}_{h}(v)(\eta, h)\right|\leq C h^{-M}\langle \eta\rangle^{M}$ for some 
$C>0,$ $M>0,$ and therefore 
$\left|\mathcal{F}_{h}(\chi_{1}u\otimes\chi_{2}v)(\xi, \eta, 
h)\right|=\mathcal{O}(h^{\infty})$ uniformly in $(\xi, \eta)\in O'_{1}\times O'_{2}$ for any 
open bounded $O'_{2}\subset\mathbb{R}^{d_{2}}$ with $\eta\in O'_{2}.$
Therefore $(x_0, \xi_0; y_0, \eta_0)\notin WF_{h}^{f}(u\otimes v).$

The proof of the second assertion in (b) is as in the $C^{\infty}$ case.
See Theorem 8.2.9, \cite{H}.

To establish (c) and (d), we first prove the following
\begin{Lem}\label{intdist}
Let $u_{1}\in
\mathcal{D}'_{h}\left(\mathbb{R}^{n}\right)$ and
$u_{2}\in\mathcal{E}'_{h}\left(\mathbb{R}^{n}\right)$ satisfy
$WF_{h}\left(u_{1}\right)\cap
WF_{h}'\left(u_{2}\right)=\emptyset.$

Then $\int 
u_{1}u_{2}=\mathcal{O}\left(h^{\infty}\right),$ where the integral is 
defined as in Proposition 7.6, \cite{GS}.
\end{Lem}

\medskip
\noindent
{\it Proof:}\;
For $u\in\mathcal{D}_{h}'(\mathbb{R}^{n})$ let
\begin{equation*}
\begin{aligned}
\Sigma_{h} & =\{\xi\in\mathbb{R}^{n}: \exists\: x\in\mathbb{R}^{n},\; (x, 
\xi)\in 
WF_{h}(u)\}\\
\Sigma_{h}^{i} & =\{\xi\in\mathbb{R}^{n}: \exists\: x\in\mathbb{R}^{n},\; 
(x, \xi)\in
WF_{h}^{i}(u)\}\\
\Sigma_{h}^{f} & =\{\xi\in\mathbb{R}^{n}: \exists\: x\in\mathbb{R}^{n},\; 
(x, \xi)\in
WF_{h}^{f}(u)\}\\
\Sigma_{h}^{x} & =\{\xi\in\mathbb{R}^{n}: (x, \xi)\in WF_{h}(u)\}\\
\end{aligned}
\end{equation*}

We have that $\Sigma_{h}^{x}(u)=\lim_{\supp \phi\to \{x\}}\Sigma_{h}(\phi 
u).$
The proof is the same as in the classical case $(h=1)$ (see \cite{H}, Section 8.1).
For every $x_{0}\in\mathbb{R}^{n}$ we can then find $\varphi\in 
C^{\infty}_{c}\left(\mathbb{R}^{n}\right)$ such that 
$\varphi\left(x_{0}\right)\ne 0$ and $\Sigma_{h}(\varphi 
u_1)\cap\left(-\Sigma_{h}(\varphi u_2)\right)=\emptyset.$
By Proposition 7.6, \cite{GS}, we have that 
\begin{equation}\label{dsjWF}
\int\varphi u_{1}\varphi u_{2} =\frac{1}{\left(2\pi 
h\right)^{n}}\int 
\widehat{\varphi 
u}_{1}\left(\frac{\xi}{h}\right)\widehat{\varphi 
u}_{2}\left(-\frac{\xi}{h}\right)d\xi.
\end{equation}

Now, since $\Sigma_{h}^{i}(\varphi u_{1})\cap\Sigma_{h}^{i}(\varphi u_{2})=\emptyset,$ for 
every 
$\xi_{0}\in
\Sigma_{h}^{i}(\varphi u_{1})$ we can find an open conic neighborhood $U_{\xi_{0}}$ of 
$\xi_{0}$ 
such that 
$\widehat{\varphi 
u_{2}}\left(-\frac{\xi}{h}\right)=\mathcal{O}(h^{\infty}\langle\xi\rangle^{-\infty})$ 
uniformly in $U_{\xi_{0}}\cap\left\{\xi: \|\xi\|\geq\frac{1}{C}\right\},$ for some $C>0.$
Since $u_{1}\in\mathcal{D}_{h}'(\mathbb{R}^{n}),$ it follows that there 
exist $N\in\mathbb{N}$ 
and $C'>0$ such that $\left|\mathcal{F}_{h}(\varphi u_{1})(\xi, h)\right|\leq C' 
h^{-N}\langle\xi\rangle^{N}$ and therefore $\widehat{\varphi 
u_{1}}\left(\frac{\xi}{h}\right) 
\widehat{\varphi 
u_{2}}\left(-\frac{\xi}{h}\right)=\mathcal{O}(h^{\infty}\langle\xi\rangle^{-\infty})$ 
uniformly in $U_{\xi_{0}}\cap\left\{\xi: \|\xi\|\geq\frac{1}{C}\right\}.$
The compactness of 
$\mathbb{S}^{n-1}$ implies that we can find finitely many such neighborhoods 
$\left(U_{l}^{1}\right)_{l=1}^{L_{1}}$ and $\left(U_{l}^{2}\right)_{l=1}^{L_{2}}$ and a 
constant
$C_{1}>0$ satisfying $\Sigma_{h}^{i}(\varphi u_{1})\subset 
\cup_{l=1}^{L_{2}}U_{l}^{2}$ and $\Sigma_{h}^{i}(\varphi 
u_{2})\subset\cup_{l=1}^{L_{1}}U_{l}^{1}$ 
and such that $\widehat{\varphi u_{1}}\left(\frac{\xi}{h}\right)\widehat{\varphi 
u_{2}}\left(-\frac{\xi}{h}\right)=\mathcal{O}(h^{\infty}\langle\xi\rangle^{-\infty})$
uniformly in $U_{l}^{j}\cap\left\{\xi: \|\xi\|\geq\frac{1}{C}\right\},$ 
$l=1, \dots, L_{j},$ $j=1, 2.$
We can further arrange to have $\left(\cup_{l=1}^{L_{1}} 
U_{l}^{1}\right)\cap 
\left(\cup_{l=1}^{L_{2}} U_{l}^{2}\right)=\emptyset.$
Lastly, we choose finitely many sets 
$\left(U_{l}\right)_{l=1}^{L_{3}}$ such that $\mathbb{S}^{n-1}\left\backslash 
\left(\cup_{k=1}^{2} 
\cup_{l}^{L_{k}}U_{l}^{k}\right)\right.\subset\cup_{l=1}^{L_{3}}U_l$ 
and a constant $C_{2}>0$ such that $\widehat{\varphi 
u_{j}}\left(\frac{\xi}{h}\right)=\mathcal{O}(h^{\infty}\langle\xi\rangle^{-\infty})$ 
uniformly in 
$\xi\in \left\{\xi: \|\xi\|>\frac{1}{C_{2}}\right\}\cap U_{l},$ $j=1, 2,$ $l=1, \dots, L_3.$ 
With $C=\min\{C_1, C_2\},$ we then have
\begin{equation*} 
\int_{\left\{\xi: \; \|\xi\|>\frac{1}{C}\right\}}
\widehat{\varphi
u}_{1}\left(\frac{\xi}{h}\right)\widehat{\varphi
u}_{2}\left(-\frac{\xi}{h}\right)d\xi=\mathcal{O}(h^{\infty}).
\end{equation*}   

The same argument applied now to 
$\Sigma_{h}^{f}(\varphi u_{j})\cap\left\{\xi\in\mathbb{R}^{n}: 
\|\xi\|\leq\frac{1}{C}\right\},$ $j=1, 2,$ gives that
\begin{equation*} 
\int_{\left\{\xi:\; \|\xi\|\leq\frac{1}{C}\right\}}
\widehat{\varphi
u}_{1}\left(\frac{\xi}{h}\right)\widehat{\varphi
u}_{2}\left(-\frac{\xi}{h}\right)d\xi=\mathcal{O}(h^{\infty})
\end{equation*} 
and therefore
\begin{equation*}
\int\varphi u_{1} \varphi u_{2}=\mathcal{O}(h^{\infty}).
\end{equation*}

Choosing a locally finite partition of unity $\sum_{j=1}^{\infty} 
\varphi_{j}^{2}=1$ with each function $\varphi_{j}$ chosen as
$\varphi$ above, we have that $\left\langle 
u_{1}, 
u_{2}\right\rangle=\sum_{j=1}^{\infty}\left\langle\varphi_{j} u_{1}, 
\varphi_{2}u_{2}\right\rangle=\mathcal{O}(h^{\infty}).$ $\hfill\Box$

We, now, turn to proving (c).
The fact that $Vv\in\mathcal{D}_{h}'(\mathbb{R}^{d_{1}})$ is proved in 
\cite[Proposition A.I.13]{ChG}.
We shall now prove that 
$\left(WF_{h}^{f}\right)'(V)(WF_{h}^{f}(v))$ is a closed set.
Let $\left(\left(x_{n},
\xi_{n}\right)\right)_{n\in\mathbb{N}}\subset\left(WF_{h}^{f}\right)'(V)(WF_{h}^{f}(v))$
converge to $(x_0, \xi_0).$
For every $n\in\mathbb{N}$ let $(y_n, \eta_n)\in WF_{h}^{f}(v)$ be such
that $(x_n, \xi_n; y_n, \eta_n)\in \left(WF_{h}^{f}\right)'(V).$
Since $WF_{h}^{f}(v)$ is compact, after passing to a subsequence, we can
assume that $(y_n, \eta_n)\to (y_0, \eta_0)\in WF_{h}^{f}(v).$
Therefore $(x_n, \xi_n; y_n, \eta_n)\to (x_0, \xi_0; y_0, \eta_0)$ and
since $\left(WF_{h}^{f}\right)'(V)$ is closed, it follows that $(x_0, 
\xi_0; y_0, \eta_0)\in
\left(WF_{h}^{f}\right)'(V).$
This implies that $(x_0, \xi_0)\in 
\left(WF_{h}^{f}\right)'(V)(WF_{h}^{f}(v))$ and therefore
$\left(WF_{h}^{f}\right)'(V)(WF_{h}^{f}(v))$ is closed.

Let, now, $(x_{0},
\xi_{0})\in\left(\left(WF_{h}^{f}\right)'(V)(WF_{h}^{f}(v))\right)^{c}$
and let $O, O'\subset\mathbb{R}^{d_{1}}$ be
open neighborhoods of $x_0$ and $\xi_0,$ respectively, such that 
\[O\times O'\subset
\left(\left(WF_{h}^{f}\right)'(V)(WF_{h}^{f}(v))\right)^{c}\]
 and $O'$ is
bounded.
Let $\chi\in C_{c}^{\infty}\left(\mathbb{R}^{d_{1}}\right)$ have support
in $O$ and let $\xi\in
O'.$
By the proof of Lemma \ref{intdist}, we have that
\begin{equation*}
\left\langle V(\cdot, \cdot\cdot), \chi(\cdot) 
e^{-\frac{i}{h}\langle\cdot,
\xi\rangle}\otimes v(\cdot\cdot) \right\rangle=\frac{1}{(2\pi
h)^{d_{1}+d_{2}}}\left\langle\hat{V}, \widehat{\chi(\cdot)
e^{-\frac{i}{h}\langle\cdot,
\xi\rangle}}\otimes\hat{v}\right\rangle=\mathcal{O}(h^{\infty})
\end{equation*}
uniformly in $\xi\in O'.$

The proof in the case of the infinite wave front set is the same.

Lastly, to prove (d), we first observe that the fact that $V\circ 
W\in\mathcal{D}_{h}'(\mathbb{R}^{d_{1}+d_{3}})$ follows as in the 
proof of \cite[Proposition A.I.13]{ChG}.
To establish (c), now, we begin by proving that 
$\left(WF_{h}^{f}\right)'\left(V\right)\circ 
\left(WF_{h}^{f}\right)'\left(W\right)$ is closed.
For that, let 
\[((x_n, \xi_n; y_n, \eta_n))_{n\in\mathbb{N}}\subset
\left(WF_{h}^{f}\right)'\left(V\right)\circ
\left(WF_{h}^{f}\right)'\left(W\right)\] converge to $(x_0, \xi_0; y_0, \eta_0).$
Let $((z_n, \zeta_n))_{n\in\mathbb{N}}\subset\mathbb{R}^{m}$ be such that $(x_n, \xi_n; z_n, 
\zeta_n)\in \left(WF_{h}^{f}\right)'(V),$ $(z_n, \zeta_n; y_n, \eta_n)\in 
\left(WF_{h}^{f}\right)'(W),$ $n\in\mathbb{N}.$
By the assumption, we can assume that, after passing to a subsequence, $(z_n, \zeta_n)\to 
(z_0, \zeta_0).$
Since then $(x_n, \xi_n; z_n, \zeta_n)\to (x_0, \xi_0; z_0, \zeta_0)$ and 
$\left(WF_{h}^{f}\right)'(V)$ is closed, it follows that $(x_0, \xi_0; z_0, \zeta_0)\in 
\left(WF_{h}^{f}\right)'(V).$
Similarly, $(y_0, \eta_0; z_0, \zeta_0)\in \left(WF_{h}^{f}\right)'(W),$ and therefore 
\[(x_0, \xi_0; y_0, \eta_0)\in \left(WF_{h}^{f}\right)'\left(V\right)\circ
\left(WF_{h}^{f}\right)'\left(W\right).\]

Let, now, $\left(x, \xi; y, \eta\right)\in 
\left(\left(WF_{h}^{f}\right)'\left(V\right)\circ 
\left(WF_{h}^{f}\right)'\left(W\right)\right)^{c}.$
Let $O_1, O'_{1}\subset\mathbb{R}^{d_{1}},$ $O_2, O'_{2}\subset\mathbb{R}^{d_{3}}$ be open 
neighborhoods of $x, \xi, y,$ and $\eta,$ respectively, such that $O_1\times O'_{1}\times 
O_2\times O'_{2}\subset \left(WF_{h}'\left(V\right)\circ WF_{h}'\left(W\right)\right)^{c}$ 
and $O'_{1}$ and $O'_{2}$ are bounded.
Let $\varphi\in C_{c}^{\infty}\left(\mathbb{R}^{d_{1}}\right),$ $\psi\in 
C_{c}^{\infty}\left(\mathbb{R}^{d_{3}}\right)$  have supports inside $O_1$ and $O_2,$ 
respectively.
Then, by the proof of Lemma \ref{intdist}, we have that
\begin{equation*}
\left\langle V(\cdot, \cdot\cdot)\otimes \psi(\cdot\cdot\cdot) 
e^{-\frac{i}{h}\left\langle \cdot\cdot\cdot, \eta\right\rangle}, 
\varphi(\cdot) e^{-\frac{i}{h}\left\langle \cdot, 
\xi\right\rangle}\otimes 
W(\cdot\cdot, \cdot\cdot\cdot)\right\rangle=\mathcal{O}\left(h^{\infty}\right) 
\end{equation*}
uniformly in $\left(\xi, \eta\right)\in O'_{1}\times O'_{2}.$
Therefore $\left(x, \xi; y, \eta\right)\notin WF_{h}^{f}(V\circ W).$ 

The proof in the infinite case is the same as in the $C^{\infty}$ case, see Theorem 7.10, 
\cite{GS}. 
$\hfill\Box$

\section{Global Semi-Classical Fourier Integral Operators}\label{shFIO}

Here we prove a characterization of global semi-classical Fourier
Integral Operators, which is the semi-classical analog of Melrose's
characterization
of Lagrangian distributions in \cite[Definition 25.1.1]{H}.

\subsection{Parametrizing Lagrangian Submanifolds}\label{sparamL}
We first review some facts from symplectic geometry relating 
non-degenerate phase functions and Lagrangian submanifolds.

Let $V\subset \mathbb{R}^{n}\times\mathbb{R}^{m},$ $m\in\mathbb{N}_{0},$
be an
open set and let
$\varphi=\varphi\left(x, \theta\right)\in C^{\infty}_{b}\left(V;
\mathbb{R}\right).$
For $m>0$, let $\varphi$ also be a phase function in the sense of
\cite{M}, Section 2.4.
If $a\in S_{n+m}\left(1\right)$, we define the
oscillatory integral
$I\left(a,\varphi\right)=\int_{\mathbb{R}^{m}}e^{\frac{i}{h}\varphi
\left(\cdot,\theta\right)}a\left(\cdot,\theta\right)d\theta$
as in
\cite{M},
Section 2.4 if $m>0$ and set $ I\left(a,\varphi\right)=
e^{\frac{i}{h}\varphi}a$ if
$m=0$.

We further let $$C_{\varphi}=\left\{\left(x,\theta\right)\in
V:\varphi'_{\theta}\left(x,\theta\right)=0\right\}$$ and
$$\Lambda_{\varphi}=\left\{\left(x,\varphi'_{x}\left(x,\theta\right)\right):
\left(x,\theta\right)\in
C_{\varphi}\right\}$$ and recall that a phase function $\varphi$ is
non-degenerate if
\begin{equation}\label{ndeg}
\varphi'_{\theta}\left(x,\theta\right)=0\text{ implies that }
\left(\varphi''_{\theta
x} \;\; \varphi''_{\theta
\theta}\right) \text{ has maximum rank at } \left(x,\theta\right).
\end{equation}
If $m=0$, it is a standard fact from symplectic geometry that
$\Lambda_{\varphi}$ is a Lagrangian submanifold of
$T^{*}\left(\mathbb{R}^{n}\right)$.
If $m>0$, (\ref{ndeg}) implies that $C_{\varphi}$ is a smooth
$n-$dimensional manifold.
Let
$j_{\varphi}:C_{\varphi}\ni\left(x,\theta\right)\mapsto\left(x,\varphi'_{x}
\left(x,\theta\right)\right)\in\Lambda_{\varphi}$.
Then, after shrinking $V$ around any fixed point
$\left(x',\theta'\right)\in
C_{\varphi}$, we can assume that $\Lambda_{\varphi}$ is a Lagrangian
submanifold of $T^{*}\mathbb{R}^{n}$ and $j_{\varphi}$ is a
diffeomorphism.
For a proof, we refer the reader to \cite{GS}, Lemmas 11.2 and 11.3.

If $\Lambda\subset T^{*}\mathbb{R}^{n}$ is a Lagrangian
submanifold such that the map $\pi_{\xi}: \Lambda\cap
U\ni\left(x,\xi\right)\mapsto\xi\in\mathbb{R}^{n}$
is a local diffeomorphism, then there exist an open set
$W\subset\mathbb{R}^{n}\backslash\left\{0\right\}$ and a function $H\in
C_{b}^{\infty}\left(W;
\mathbb{R}\right)$ satisfying
\begin{equation}\label{H}
\Lambda\cap U=\left\{\left(H'\left(\xi\right),\xi\right):\xi\in W\right\}.
\end{equation}
For a proof, see \cite{GS}, Section 9.

If $\Lambda\subset T^{*}\mathbb{R}^{n}$ is any Lagrangian submanifold and
$\gamma\in\Lambda$, then there exists an open set $U\subset
T^{*}\mathbb{R}^{n}$ and a
non-degenerate phase function $\varphi\in C^{\infty}\left(V\right),$
$V\subset\mathbb{R}^{n+m}$
open, $m\in\mathbb{N}_{0}$ such that
\begin{equation}\label{lphi}
\Lambda\cap U=\Lambda_{\varphi}.
\end{equation}
We include the proof of this well-known result here for completeness and
to introduce some
notation.
Let $\mu=T_{\gamma}\Lambda$ be identified in a natural way with a subspace
of
$T^{*}\mathbb{R}^{n}.$
By Lemma 9.5, \cite{GS} we have that after a linear change of coordinates
we may
assume that
\begin{equation}\label{tplane}
\mu=\left\{\left(0, x''; \xi', B x''\right)\right\},
\end{equation}
for a splitting of the coordinates $x=\left(x', x''\right)$ and
$\xi=\left(\xi', \xi''\right),$
where $x'=\left(x_1,
\dots, x_k\right),$ $k=0, \dots, n,$ and $B$ is a real symmetric matrix.
This implies that the differential of the projection
$\pi:\Lambda\rightarrow \left(x'', \xi'\right)$ is
bijective at $\gamma$ and therefore this map is a local diffeomorphism
from a neighborhood
of $\gamma$ to the $\left(x'', \xi'\right)-$space.
Therefore there exists a function $S\in C^{\infty}\left(\mathbb{R}^{n};
\mathbb{R}\right)$ and an open
neighborhood $U\subset T^{*}\mathbb{R}^{n}$ of $\gamma$ such that
$\Lambda\cap
U=\left\{\left(\frac{\partial S}{\partial \xi'}, x''; \xi',
-\frac{\partial S}{\partial
x''}\right)\right\}\cap U.$
From this it easily follows that $\varphi\left(x,
\xi'\right)=\left\langle
x',
\xi'\right\rangle-S\left(x'', \xi'\right)$ is a
non-degenerate phase function such that $\Lambda\cap
U=\Lambda_{\varphi}\cap U.$

\subsection{Semi-Classical Fourier Intergal Operators}

We are now ready to make the following definition
\begin{Def}\label{dfio}
Let $M$ be a smooth $k$-dimensional manifold and let $\Lambda\subset
T^{*}M$ be a smooth closed
Lagrangian submanifold with respect to the canonical symplectic 
structure on $T^{*}M.$
Let $r\in\mathbb{R}.$
Then the space $I^{r}_{h}\left(M, \Lambda\right)$ of semi-classical
Fourier integral
distributions of order $r$ associated to $\Lambda$ is defined as the set
of all $u\in\mathcal{D}'_{h}\left(M\right)$ 
such
that
\begin{equation}\label{defgfio}
\left(\prod_{j=0}^{N}
A_{j}\right)\left(u\right)=\mathcal{O}_{L^{2}\left(M\right)}\left(h^{N-r-\frac{k}{4}}\right),
h\to 0,
\end{equation}
for all $N\in\mathbb{N}_{0}$ and for all $A_{j}\in \Psi_{h}^{0}\left(1,
X\right),$ $j=0, \dots, N-1,$ with
compactly
supported symbols and principal symbols vanishing on $\Lambda$, and any $ 
A_N \in 
\Psi_h^{0} ( 1 , X ) $ with a compactly supported symbol.

A continuous linear operator
$C_{c}^{\infty}\left(M_1\right)\rightarrow\mathcal{D}_{h}'\left(M_2\right),$ 
where $M_1, M_2$ are smooth manifolds, 
whose Schwartz kernel is an element of 
$I_{h}^{r}(M_1\times M_2, \Lambda)$ for some 
Lagrangian submanifold $\Lambda\subset T^{*}M_1\times T^{*}M_2$ and some $r\in\mathbb{R}$ 
will be called a global semi-classical Fourier integral
operator of order $r$ associated to $\Lambda.$
We denote the space of these operators by 
$\mathcal{I}_{h}^{r}(M_1\times M_2, \Lambda).$
\end{Def}

\medskip
\noindent
{\bf Remark:} The exotic looking numerology for the 
 order needs to be explained. We follow the same convention as that in 
classical case and require that pseudodifferential operators with compactly 
supported 
symbols in $ S^0( 1) $ have kernels in $ I^{0}_{h}({\mathbb R}^{2n}, 
N^* \Delta ) $, where $ \Delta $ is the diagonal in $ {\mathbb R}^{2n}$.
Explicitly, suppose that 
\[ K( x , y ) = \frac{1}{ ( 2 \pi h)^n } \int e^{ \frac{i}{h} 
\langle x - y , \xi \rangle } c ( x , \xi ) d\xi \,, \ \ c \in S^0( 1) 
\cap C^\infty_{\rm{c}} ( {\mathbb R}^{2n} )\,.\]
Then 
\[ \| K \|_{ L^2 ( {\mathbb R}^{2n } ) } = 
\frac{1}{(2\pi h)^{\frac{n}{2}}}\| {\mathcal F}_{h}^{y} K \|_{L^2 ( 
{\mathbb R}^{2n } )} = 
\frac{1}{(2 \pi h)^{\frac{n}{2} }} \| c \|_{L^2 ( {\mathbb R}^{2n } )} 
\,,\]
where $ {\mathcal F}_{h}^{y} $ is the semi-classical Fourier transform in 
the
$ y $ variable, 
which is consistent with \eqref{defgfio} with $ N=0 $ and the order 
$ r = 0 $ ($ k = 2n $ here).

We, now, have the following semi-classical analog of 
\cite[Lemma 25.1.2]{H}, vol. IV.
\begin{Lem}\label{gfio}
If $u\in I_{h}^{r}(M, \Lambda),$ then $Au\in I_{h}^{r}(M, \Lambda)$ for 
every $A\in\Psi_{h}^{0}(1, M)$ or $A\in \Psi_{h}^{0, k}(T^{*}M),$ 
$k\in\mathbb{R}$ with compactly supported symbol. 

If $u\in\mathcal{D}_{h}'(M)$ is such that for every $(x_{0},
\xi_{0})\in\Lambda$ there exists $A\in\Psi_{h}^{0}(1, M)$ elliptic at 
$(x_{0}, \xi_{0})$ with compactly 
supported symbol and $Au\in I_{h}^{r}(M, \Lambda),$ then $u\in 
I_{h}^{r}(M, \Lambda).$
The same conclusion holds if $A\in\Psi_{h}^{0, k}(T^{*}M),$ 
$k\in\mathbb{R}.$
\end{Lem}

\medskip
\noindent
{\it Proof:}\;
To prove the first statement, let $u\in I_{h}^{r}(M, 
\Lambda),$ 
let 
$A\in\Psi_{h}^{0}(1, M)$ have a compactly supported symbol, and let 
$A_{j}\in 
\Psi_{h}^{0}(1, M),$ $j=1, \dots, N,$ $N\in\mathbb{N}$, also have compactly supported 
symbols and principal symbols vanishing on $\Lambda.$
Then
\begin{equation*}
\left(\prod_{j=1}^{N}A_{j}\right)(Au)=\left(\prod_{j=1}^{N-1}A_{j}\right)[A_{N}, 
A]u+\left(\prod_{j=1}^{N-1}A_{j}\right)AA_{N}u.
\end{equation*}
Here $[A_{N}, A]\in\Psi_{h}^{-1}(1, M)$ has a compactly supported symbol and 
therefore \[\left\|[A_{N}, 
A]u\right\|_{L^{2}(M)}=\mathcal{O}\left(h^{1-r-\frac{k}{4}}\right).\]
By the choice of $A$ and $A_{N}$ we further have that \[\|A A_N 
u\|_{L^{2}(M)}=\mathcal{O}\left(h^{1-r-\frac{k}{4}}\right).\]
Thus, if follows by induction with respect to $N$ that 
\[\left(\prod_{j=1}^{N}A_{j}\right)(Au)=\mathcal{O}_{L^{2}(M)}\left(h^{N-r-\frac{k}{4}}\right),\; h\to 0.\]
Therefore $Au\in I_{h}^{r}(M, \Lambda).$

To prove the converse, let $B\in \Psi_{h}^{0}(1, M)$ have a compactly 
supported symbol and satisfy $(x_{0}, 
\xi_{0})\notin WF_{h}(BA-I).$
Then $(x_{0}, \xi_{0})\notin WF_{h}(BAu-u).$ 
From the first part of the proof, we have that $BAu\in I_{h}^{r}(M, 
\Lambda).$
Let, now, $P\in\Psi_{h}^{0}(1, M)$ have symbol supported in a sufficiently 
small neighborhood of $(x_0, \xi_0)\in\Lambda$ so that 
$PBAu-Pu=\mathcal{O}(h^{\infty})$ in $C^{\infty}(M).$
Since again $PBAu\in I_{h}^{r}(M, \Lambda),$ we have that 
$\left(\prod_{j=1}^{N}A_{j}\right)(Pu)=\mathcal{O}_{L^{2}}\left(h^{N-r-\frac{k}{4}}\right),$ 
$h\to 0,$ for any set of operators $\left( A_{j}\right)_{j=1}^{N},$ $N\in\mathbb{N},$ as in 
(\ref{defgfio}).
Thus $Pu\in I_{h}^{r}(M, \Lambda)$ for every $P\in\Psi_{h}^{0}(1, M)$ with 
symbol supported in a sufficiently small neighborhood of any point 
$(x_{0}, \xi_{0})\in \Lambda.$
The compactness of the supports of the operators $\left( A_{j}\right)_{j=1}^{N},$ $N\in\mathbb{N},$ now allows us to find $P_{j}\in\Psi_{h}^{0}(1, M),$ $j=1, \dots, J,$ 
$J\in\mathbb{N},$ such that $P_{j}u\in I^{r}_{h}(M, \Lambda),$ $j=1, \dots, 
J,$ and $\sum_{j=1}^{J}\sigma(P_j)=1$ on $\cup_{k=1}^{N}\supp 
\sigma(A_{k}).$
Using the calculus of semi-classical pseudodifferential operators, we further obtain 
\begin{equation*}
\mathcal{O}_{L^{2}(M)}\left(h^{N-r-\frac{k}{4}}\right)=
\left(\prod_{k=1}^{N}A_{k}\right)\left(\sum_{j=1}^{J}P_{j}\right)u
=\left(\prod_{k=1}^{N}A_{k}\right)u
+\mathcal{O}_{L^{2}(M)}(h^{\infty}),
\end{equation*}
which completes the proof. 

The proof in the case of an operator $A\in\Psi_{h}^{0, k}(T^{*}M),$
$k\in\mathbb{R}$ is analogous. $\hfill\Box$

We shall now characterize semi-classical Fourier integral distributions 
microlocally.
We have the following

\begin{Th}\label{lfio}
Let $\Lambda \subset T^{*}\mathbb{R}^{n}$ be a Lagrangian submanifold and
let $\gamma\in \Lambda.$
Let $\varphi$ be a non-degenerate phase function in an open set 
$V\subset\mathbb{R}^{n+m},$
$m\in\mathbb{N}_{0},$ such that $\Lambda=\Lambda_{\varphi}$ in a 
neighborhood
of $\gamma.$
If $a\in S_{n+m}^{r+\frac{m}{2}+\frac{n}{4}}\left(1\right)$ is such that 
$\supp 
a\Subset V$, then
$I\left(a,\varphi\right)\in I_{h}^{r}\left(\mathbb{R}^{n}, 
\Lambda\right).$

Conversely, if $u\in I_{h}^{r}\left(\mathbb{R}^{n}, \Lambda\right)$ 
microlocally 
near $\gamma$, then for every non-degenerate phase function $\varphi$ in 
an open set
$V\subset\mathbb{R}^{n+m},$ $m\in\mathbb{N}_{0},$ such that
$\Lambda=\Lambda_{\varphi}$ near $\gamma$, there exists
$a\in S_{n+m}^{r+\frac{m}{2}+\frac{n}{4}}\left(1\right)$ with $\supp 
a\Subset V$
such that $u=I\left(a,\varphi\right)$ microlocally near $\gamma$.
\end{Th}

\medskip
\noindent
{\it Proof:}\;
Let $\gamma$ have canonical coordinates $\left(x_{0}, \xi_{0}\right)$ and 
let us 
first assume that 
$\Lambda$ is transverse to the section $\xi=\xi_{0}$ at $\gamma.$
Then there exists an open neighborhood $U\subset T^{*}\mathbb{R}^{n}$ of 
$\gamma$ such that $\pi_{\xi}: \Lambda\cap U\ni \left(x, \xi\right)\mapsto 
\xi\in 
\mathbb{R}^{n}$ in 
canonical coordinates is a local diffeomorphism.
Let $H\in C^{\infty}_{b}\left(\mathbb{R}^{n};\mathbb{R}\right)$ be chosen 
such that, 
perhaps after 
adjusting $V,$ 
$\Lambda_{\varphi}=\left\{\left(H'\left(\xi\right),\xi\right):\xi\in 
W\right\}$ for 
some 
bounded open 
set $W\subset\mathbb{R}^{n}.$ 
For  $\xi\in \mathbb{R}^{n}$ consider
\begin{displaymath}
\widehat{I\left(a,\varphi\right)}\left(\xi\right)=\int \int 
e^{\frac{i}{h}\left(\varphi\left(x,\theta\right)-\left\langle
x,\xi\right\rangle\right)} a\left(x,\theta\right)d\theta dx.
\end{displaymath}

For $\xi\not\in W$ integration by parts in $\left(x,\theta\right)$ gives
\begin{equation}\label{onwc}
\widehat{I\left(a,\varphi\right)}=\mathcal{O}\left(h^{\infty}\right) 
\text{ in } 
C_{c}^{\infty}\left(W^{c}\right).
\end{equation}

Let, now, $\bar{\xi}\in W$.
Then the function $\Phi\left(x,\theta; 
\bar{\xi}\right)=\varphi\left(x,\theta\right)-\left\langle 
x, \bar{\xi} \right\rangle$ has a critical point at 
$\left(\bar{x}\left(\bar{\xi}\right),\bar{\theta}\left(\bar{\xi}\right)\right)$ 
which is the inverse image 
in 
$C_{\varphi}$ under $j_{\varphi}$ of the point 
$\left(H'\left(\bar{\xi}\right), \bar{\xi}\right)$.
Using integration by parts again, we obtain that, up to a term which is 
$\mathcal{O}\left(h^{\infty}\right)$ in $C_{c}^{\infty}\left(W\right)$,
\begin{equation}\label{hft}
\widehat{I\left(a,\varphi\right)}\left(\bar{\xi}\right)\equiv\int \int 
e^{\frac{i}{h}\left(\varphi\left(x,\theta\right)-\left\langle
x, \bar{\xi} \right\rangle\right)} 
a\left(x,\theta\right)\chi\left(x-\bar{x}\left(\bar{\xi}\right), 
\theta-\bar{\theta}\left(\bar{\xi}\right)\right)d\theta 
dx,
\end{equation}
where $\chi\in C_{c}^{\infty}\left(\mathbb{R}^{n+m}\right)$ is equal to 1 
on a 
neighborhood of 0.

To prove that the critical point is non-degenerate, let $v$ be in the 
kernel of 
\begin{displaymath}
\Phi''_{x\theta}\left(\bar{x}\left(\bar{\xi}\right),\bar{\theta}\left(\bar{\xi}\right); 
\xi_{0}\right)=\begin{bmatrix} 
\varphi''_{xx}\left(\bar{x}\left(\bar{\xi}\right),\bar{\theta}\left(\bar{\xi}\right)\right)
&\varphi''_{x\theta}\left(\bar{x}\left(\bar{\xi}\right), 
\bar{\theta}\left(\bar{\xi}\right)\right)\\ 
\varphi''_{\theta x}\left(\bar{x}\left(\bar{\xi}\right), 
\bar{\theta}\left(\bar{\xi}\right)\right) & 
\varphi''_{\theta 
\theta}\left(\bar{x}\left(\bar{\xi}\right), 
\bar{\theta}\left(\bar{\xi}\right)\right)
\end{bmatrix}.
\end{displaymath}
Then $$v\in \ker \left(\varphi''_{x 
\theta}\left(\bar{x}\left(\bar{\xi}\right), 
\bar{\theta}\left(\bar{\xi}\right)\right) \;\; 
\varphi''_{\theta
\theta}\left(\bar{x}\left(\bar{\xi}\right), 
\bar{\theta}\left(\bar{\xi}\right)\right)\right)$$ and therefore $$v\in 
T_{\left(\bar{x}\left(\bar{\xi}\right), 
\bar{\theta}\left(\bar{\xi}\right)\right)}C_{\varphi}.$$
We also have that $$v\in\ker \left(\varphi''_{xx}\left(\bar{x}\left(\bar{\xi}\right), 
\bar{\theta}\left(\bar{\xi}\right)\right) \;\; 
\varphi''_{x\theta}\left(\bar{x}\left(\bar{\xi}\right), 
\bar{\theta}\left(\bar{\xi}\right)\right)\right)$$ and 
since 
$j_{\varphi}$ 
and $\pi_{\xi}$ are diffeomorphisms, it follows that $v=0.$
Hence 
\begin{equation}\label{ns}
\text{ the matrix } \Phi''_{x\theta}\left(\bar{x}\left(\bar{\xi}\right), 
\bar{\theta}\left(\bar{\xi}\right); 
\bar{\xi}\right) 
\text{ is non-singular.}
\end{equation}
We can therefore apply the method of stationary phase to the integral 
(\ref{hft}) and obtain
\begin{equation}\label{sph}
\widehat{I\left(a, \varphi\right)}\left(\bar{\xi}\right)\sim
e^{\frac{i}{h}\Phi\left(\bar{x}\left(\bar{\xi}\right), 
\bar{\theta}\left(\bar{\xi}\right);\bar{\xi}\right)}\sum_{k=0}^{\infty}h^{k+\frac{n}{2}+\frac{m}{2}}\left(A_{2k}\left(D_{x,\theta}\right)a\right)\left(\bar{x}\left(\bar{\xi}\right), 
\bar{\theta}\left(\bar{\xi}\right)\right), 
\end{equation}
where $A_{2k}\left(D\right)$ are differential operators of orders $\leq 
2k$, 
respectively.

The Implicit Function Theorem and (\ref{ns}) now imply that, 
perhaps after 
shrinking $W$ around $\bar{\xi}$, $\bar{x}, \bar{\theta}\in 
C^{\infty}\left(W\right)$.
We further adjust $W$ so that $\bar{x}, \bar{\theta}\in 
C^{\infty}_{b}\left(W\right).$
Thus 
\[\Phi'_{\xi}\left(\bar{x}\left(\xi\right),\bar{\theta}\left(\xi\right);\xi\right)=-H'\left(\xi\right),\; 
\xi\in W\] and therefore, by adding a constant to $H$ if necessary, we can assume 
that
$\Phi\left(\bar{x}\left(\xi\right), 
\bar{\theta}\left(\xi\right);\xi\right)=-H\left(\xi\right)$ for $\xi\in 
W$.
We also have that for every $k$, 
$A_{2k}\left(D_{x,\theta}\right)\left(a\right)\left(\bar{x}\left(\cdot\right), 
\bar{\theta}\left(\cdot\right)\right)\in 
S_{n}^{r+\frac{m}{2}+\frac{n}{4}}\left(1\right)$.
Thus, with $A\in S_{n}^{r-\frac{n}{4}}\left(1\right)$, $A\sim 
\sum_{k=0}^{\infty}h^{k+\frac{n}{2}+\frac{m}{2}}\left(A_{2k}\left(D_{x,\theta}\right)a\right)\left(\bar{x}\left(\cdot\right), 
\bar{\theta}\left(\cdot\right)\right),$ 
$A=\mathcal{O}\left(h^{\infty}\right)$ in 
$S_{n}^{r-\frac{n}{4}}\left(1\right)$ outside $W$, we 
obtain, from 
(\ref{onwc}), (\ref{hft}), and 
(\ref{sph}), that
\begin{displaymath}
\widehat{I\left(a, 
\varphi\right)}\left(\xi\right)=e^{-\frac{iH\left(\xi\right)}{h}}A\left(\xi\right).
\end{displaymath}

Now, the ideal of smooth functions vanishing on $\Lambda_{\varphi}$ is 
generated by the symbols 
$a_{j}\left(x,\xi\right)=x_{j}-H'_{\xi_{j}}\left(\xi\right)$, 
$j=1,\dots, n.$
Since $I\left(a, \varphi\right)$ has compact support, by adjusting $V$ 
without 
changing
$I\left(a, \varphi\right)$, we can assume that $\Lambda_{\varphi}$ is 
compact and we 
can choose
$\chi\in C_{c}^{\infty}\left(T^{*}\mathbb{R}^{n}\right)$ equal to 1 on a 
neighborhood
of $\Lambda_{\varphi}$.
Then $\tilde{a}_{j}=\chi a_{j}\in S_{2n}\left(1\right)$, $j=1, \dots, n$ 
vanish on
$\Lambda_{\varphi}$.
By the calculus of pseudodifferential operators, we have that
$Op_{h}\left(\tilde{a}_{j}\right)I\left(a, 
\varphi\right)=Op_{h}\left(a_{j}\right)I\left(a, 
\varphi\right)+E_{j}I\left(a, 
\varphi\right)$, where
$E_{j}I\left(a, \varphi\right)=\mathcal{O}\left(h^{\infty}\right)$, $h\to 
0$ in 
$C_{c}\left(\mathbb{R}^{n}\right)$.
\begin{equation*}
\begin{aligned}
& \mathcal{O}\left(h^{\infty}\right)+ 
\left\|\left(Op_{h}\left(\tilde{a}_{j}\right)\right)^{\alpha}\left(I\left(a, 
\varphi\right)\right)\right\|_{L^{2}\left(\mathbb{R}^{n}\right)}=
\left\|\left(x-H'\left(hD\right)\right)^{\alpha}\left(I\left(a, 
\varphi\right)\right)\right\|_{L^{2}\left(\mathbb{R}^{n}\right)}\\
 & \quad = \frac{1}{\left(2\pi 
h\right)^{\frac{n}{2}}}\left\|\left(-hD-H'\right)^{\alpha}\widehat{I\left(a, 
\varphi\right)}\right\|_{L^{2}\left(\mathbb{R}^{n}\right)}=\mathcal{O}\left(h^{|\alpha|-r-\frac{n}{4}}\right), 
\alpha\in \mathbb{N}^{n}, h\to 0.\\
\end{aligned}
\end{equation*}
Thus $I\left(a, \varphi\right)\in I_{h}^{r}\left(\mathbb{R}^{n}, 
\Lambda\right).$ 

We remark here that the same argument will allow us in similar situations
to use in condition (\ref{defgfio}) symbols, which do not 
belong to the 
class $S\left(1\right)$
and below we will do so without repeating this argument.

We now turn to proving the converse.
Let $U, H, W,$ and  $V$ be further chosen so that (\ref{H})
and (\ref{lphi}) hold and $W$ is bounded.
Extend $H$ to a function in 
$C^{\infty}_{b}\left(\mathbb{R}^{n};\mathbb{R}\right)$.
Let $P\in \Psi^{0}_{h}\left(1, \mathbb{R}^{n}\right)$ satisfy 
(\ref{P}) and set 
$\tilde{u}=Pu$.
The symbols $a_{j}\left(x,\xi\right)=x_{j}-H'_{\xi_{j}}\left(\xi\right)$, 
$j=1, \dots, n$ vanish 
on $\Lambda_{\varphi}\cap U$, and therefore we obtain from 
(\ref{defgfio}) 
that
\begin{displaymath}
\left\|\left(x-H'\left(hD\right)\right)^{\alpha}\left(\tilde{u}\right)\right\|_{L^{2}\left(\mathbb{R}^{n}\right)}=\mathcal{O}\left(h^{|\alpha|-r-\frac{n}{4}}\right), 
\alpha\in \mathbb{N}^{n}, h\to 0
\end{displaymath}
and hence, after taking the Fourier transform,
\begin{equation}\label{comp}
\left\|\left(-hD-H'\right)^{\alpha}\hat{\tilde{u}}\right\|_{L^{2}\left(\mathbb{R}^{n}\right)}=\mathcal{O}\left(h^{|\alpha|+\frac{n}{4}-r}\right), 
\alpha\in \mathbb{N}^{n}, h\to 0.
\end{equation}
Substituting 
$\hat{\tilde{u}}\left(\xi\right)=e^{-\frac{iH\left(\xi\right)}{h}}v\left(\xi\right)$ 
in 
(\ref{comp}), we obtain
\begin{equation}
\begin{aligned}
\left\|\left(hD\right)^{\alpha}v\right\|_{L^{2}\left(\mathbb{R}^{n}\right)}&=\mathcal{O}\left(h^{|\alpha|+\frac{n}{4}-r}\right), 
\alpha\in 
\mathbb{N}^{n}, h\to 0,\\
\left\|D^{\alpha} v
\right\|_{L^{2}\left(\mathbb{R}^{n}\right)}&=\mathcal{O}\left(h^{\frac{n}{4}-r}\right), 
\alpha\in 
\mathbb{N}^{n}.
\end{aligned}
\end{equation}
As in the previous case, this implies that
\begin{displaymath}
v\in S_{n}^{r-\frac{n}{4}}\left(1\right).
\end{displaymath}
Let $\Phi\left(x,\theta;\xi\right)=\varphi 
\left(x,\theta\right)-\left\langle x,\xi\right\rangle$, $\xi\in 
W$, $\left(x,\theta\right)\in V$.
Choose $\bar{\xi}\in W$, and let $\left(\bar{x}\left(\bar{\xi}\right), 
\bar{\theta}\left(\bar{\xi}\right)\right)\in 
C_{\varphi}$ be the critical point of $\Phi\left(\cdot, \cdot\cdot; 
\bar{\xi}\right)$.
Let $M\Subset V$ be a neighborhood of 
$\left(\bar{x}\left(\bar{\xi}\right), 
\bar{\theta}\left(\bar{\xi}\right)\right)$ 
such that $\text{sgn } \Phi''$ is constant on $M$ 
and let $\psi\in C^{\infty}\left(M; \mathbb{R}^{n}\right)$ be such that 
$\psi\left(x,\theta\right)=\varphi'_{x}\left(x,\theta\right)$ on 
$C_{\varphi}$.
Define 
$$a_{0}\left(x,\theta\right)=\frac{1}{\left(2\pi 
h\right)^{\frac{n}{4}+\frac{m}{2}}}\left(e^{-\frac{i\pi}{4}\text{sgn} 
\Phi''_{x\theta}}|\text{det}\Phi''_{x\theta}|^{\frac{1}{2}}v\right)\circ 
\psi\left(x,\theta\right) \text{ for } \left(x,\theta\right)\in M.$$
Then, by the first part of the proof, we have that 
\begin{displaymath}
e^{-\frac{iH}{h}}v-\widehat{I\left(a_{0},\varphi\right)}=\mathcal{O}\left(h\right), 
h\to 0, 
\text{ in } S_{n}^{r-\frac{n}{4}}\left(1\right).
\end{displaymath} 
Iterating this process, we obtain a sequence of symbols $a_{l}\in 
S_{n}^{r+\frac{m}{2}+\frac{n}{4}}\left(1\right)$ such that, if we denote 
$U_{s}= 
I\left(\sum_{l=0}^{s} 
h^{l}a_{l}, 
\varphi\right),$ $s\in\mathbb{N}_{0}$ we have that
\begin{displaymath}
e^{-\frac{iH}{h}}v-\hat{U}_{s}=\mathcal{O}\left(h^{s+1}\right), h\to 0, 
\text{ in } 
S_{n}^{r-\frac{n}{4}}\left(1\right).
\end{displaymath}
Therefore, if we choose an asymptotic sum $a\in 
S_{n+m}^{r+\frac{m}{2}+\frac{n}{4}}\left(1\right)$ of 
$\sum_{k=0}^{\infty} h^{k}a_{k}$, we obtain $\tilde{u}=I\left(a, 
\varphi\right)$ 
microlocally near $\gamma$.

Let, now, $\Lambda\subset T^{*}\mathbb{R}^{n}$ be any Lagrangian 
submanifold and assume 
that the coordinates have been chosen in such a way that 
$\mu=T^{*}_{\gamma}\Lambda$ has the form (\ref{tplane}). 
Choose a real symmetric matrix $A_{\Lambda}=\begin{bmatrix}
0_{k\times k} & 0\\ 0 & D_{\Lambda \; 
\left(n-k\right)\times\left(n-k\right)} \end{bmatrix}$
such that 
\begin{equation}\label{det}
\det\left(B+D_{\Lambda}\right)\ne 0.
\end{equation}
Let $\tilde{\Lambda}=\left\{\left(x, \xi+A_{\Lambda}x\right): \left(x, 
\xi\right)\in\Lambda\right\}$ and 
let $\tilde{\gamma}=\left(x_{0}, 
\xi_{0}+A_{\Lambda}x_{0}\right)=\left(x_{0}, \eta_{0}\right),$ where 
$\left(x_{0}, \xi_{0}\right)$ are 
the coordinates of $\gamma.$
Then, if $\varphi\in C^{\infty}\left(V; \mathbb{R}\right),$ $V\subset 
\mathbb{R}^{n+m},$
$m\in\mathbb{N}_{0},$ is a non-degenerate phase function which 
parameterizes $\Lambda$ near 
$\gamma,$ it is clear that $\tilde{\varphi}\left(x, 
\theta\right)=\frac{1}{2}\left\langle A_{\Lambda} x, x \right\rangle + 
\varphi\left(x, \theta\right)$ 
is a non-degenerate phase 
function which parameterizes $\tilde{\Lambda}$ near $\tilde{\gamma}.$

Let $\tilde{\mu}=T_{\tilde{\gamma}}\tilde{\Lambda}.$
It is easy to see that $\tilde{\mu}=\left\{\left(0, x''; \xi', 
\left(B+D_{\Lambda}\right)x''\right)\right\}$ and it then 
follows from (\ref{det}) that $\tilde{\Lambda}$ is transverse 
to the 
constant section 
$\eta=\eta_{0}$ at $\tilde{\gamma}.$ 

Let $u\in I_{h}^{r}\left(\mathbb{R}^{n}, \Lambda\right)$ microlocally near 
$\gamma$ 
and let $A_k\in\Psi_{h}^{0}\left(1, X\right),$ $k=1, \dots, N,$ be such 
that 
$\sigma_{0}\left(A_k\right)|_{\Lambda}=0.$
From (\ref{defgfio}) we have that
\begin{equation*}
\left(\prod_{k=1}^{N} e^{\frac{i}{h}\left\langle A_{\Lambda}\cdot, \cdot 
\right\rangle}A_{k} 
e^{-\frac{i}{h}\left\langle A_{\Lambda}\cdot\cdot, \cdot\cdot 
\right\rangle} \right) e^{\frac{i}{h}\left\langle A_{\Lambda} \cdot\cdot, 
\cdot\cdot 
\right\rangle 
}u=\mathcal{O}_{L^{2}(\mathbb{R}^{n})}\left(h^{N-r}\right), h\to 0.
\end{equation*}
Let $K\left(z, t\right)=\frac{1}{\left(2\pi 
h\right)^{n}}e^{-\frac{i}{h}\left\langle A_{\Lambda}z, 
z\right\rangle}\int 
e^{\frac{i}{h}\left\langle z-t, \xi\right\rangle}a\left(z, \xi\right)d\xi 
\; 
e^{\frac{i}{h}\left\langle A_{\Lambda}t, 
t\right\rangle}$ and consider 
\begin{equation*}
\begin{aligned}
b\left(w, \tau\right) & =\int e^{\frac{i}{h}y\cdot\tau}K\left(w, 
w+y\right)dy\\
 & =\frac{1}{\left(2\pi h\right)^{n}}\int\int 
e^{\frac{i}{h}\left(\left\langle y,\tau-\xi 
\right\rangle-\left\langle 
A_{\Lambda}w, w\right\rangle + \left\langle A_{\Lambda}\left(w+y\right), 
\left(w+y\right)\right\rangle\right)}a\left(w, \xi\right) 
d\xi dy.
\end{aligned}
\end{equation*}
We apply the method stationary phase to the above integral.
The phase has a critical point at $\left(y_0, \xi_0\right)=\left(0, 
\tau+A_{\Lambda}w\right)$ 
and the Hessian of the phase 
at the critical point is $\begin{bmatrix} A_{\Lambda} & -I\\-I & 
0\end{bmatrix},$ which has 
determinant 1 and signature 0.
Therefore $b\sim \sum_{k=0}^{\infty}h^{k}b_{k}$ with $b_k\in 
S_{2n}\left(1\right)$ 
and 
\begin{equation}\label{ba}
b_{0}\left(w, \tau\right)=a\left(w, \tau+A_{\Lambda}w\right).
\end{equation}
This implies that $b\in S_{2n}\left(1\right)$ and $B_k= 
e^{\frac{i}{h}\left\langle 
A_{\Lambda}\cdot, \cdot 
\right\rangle}A_{k} e^{-\frac{i}{h}\left\langle A_{\Lambda}\cdot\cdot, 
\cdot\cdot 
\right\rangle}\in\Psi_{h}^{0}\left(1, X\right).$
From (\ref{ba}) we then have that 
$\sigma_{0}\left(B_k\right)|_{\tilde{\Lambda}}=0,$ 
$k=1, \dots, N.$
We can now apply the first part of the proof of this theorem and we have 
that 
$e^{\frac{i}{h}\left\langle A_{\Lambda}\cdot, \cdot 
\right\rangle}u=I\left(a, \tilde{\varphi}\right),$ $a\in 
S_{n+m}^{r+\frac{m}{2}+\frac{n}{4}}\left(1\right).$
Therefore $u=I\left(a, \varphi\right).$

The converse follows from reversing this argument. \hfill $\Box$

\medskip
\noindent {\bf Remark.}  Let $u\in I^{r}_{h}\left(M, \Lambda\right).$ 
Then Theorem \ref{lfio} and Lemma \ref{gfio} implies that for any 
$P\in\Psi^{0}_{h}\left(1, 
M\right)$ with 
a 
compactly supported symbol $Pu$ is given by a finite sum of 
oscillatory integrals of the form $h^{-r}I\left(a, \varphi\right),$ where 
$a\in S\left(1\right)$ 
and $\varphi$ is 
a non-degenerate phase function such that $\Lambda=\Lambda_{\varphi}$ near 
a point in 
$\Lambda.$\hfill$\Box$

Following this remark, we see that after taking a locally finite partition of unity 
$\left(a_j\right)_{j=1}^{\infty}\subset C_{c}^{\infty}(T^{*}M)$ such that 
$\sum_{j=1}^{\infty} a_j=1$ on $\Lambda$ and applying an integration by parts argument, as 
in \cite{GS}, Chapter 7, we have

\begin{Lem}\label{wfhfio}
If $u\in I_{h}^{r}(M, \Lambda),$ then $WF_{h}^{f}(u)\subset\Lambda.$
\end{Lem}

\subsection{Generalization of Egorov's Theorem}

We now prove the following generalization of Egorov's Theorem to manifolds of unequal dimensions:

\begin{Lem}\label{interwine}
Let $X_{j}$, $j=1, 2,$ be smooth manifolds.
Let $\sigma_{j}$ be the canonical symplectic form on $T^{*}X_{j}$, and 
$\pi_{j}:
T^{*}X_{1}\times T^{*}X_{2}\rightarrow T^{*}X_{j}$ the projection onto
the $j$-~th factor.
Let $$\Lambda\subset T^{*}X_{1} \times T^{*}X_{2}$$ be a Lagrangian 
submanifold of $$\left(T^{*}X_{1} \times T^{*}X_{2}, 
\pi_{1}^{*}\sigma_{1}+\pi_{2}^{*}\sigma_{2}\right)$$ such that 
$\pi_2|_{\Lambda}$ is an immersion.
Let $F\in \mathcal{I}_{h}^{r}\left(X_1\times X_2, \Lambda\right),$ 
$r\in\mathbb{R},$ have a non-vanishing principal symbol at $ ( \rho_1 , 
\rho_2 ) \in \Lambda  $. 

Then for every $A\in \Psi_{h}^{0}\left(1, X_{1}\right)$ with symbol 
supported in a 
sufficiently small neighborhood of $ \rho_{1} $ there 
exists 
$B\in 
\Psi_{h}^{0}\left(1, X_{2}\right)$ with symbol supported in a sufficiently 
small 
neighborhood of $\rho_{2}$ 
such that $$AF\equiv FB \text{ near } \left(\rho_{1}, \rho_{2}\right)$$ 
and 
\[i^{n_{2}}\left(\pi_{2}|_{\Lambda}\right)^{*}\sigma_{0}\left(B\right)=i^{n_{1}}\left(\pi_{1}|_{\Lambda}\right)^{*}\sigma_{0}\left(A\right).\]
\end{Lem}

\medskip

\noindent
{\bf Remark:} Strictly speaking we have not defined a symbol of a 
Fourier integral operator given in Definition \ref{dfio}. However,
the proof of Theorem \ref{lfio} shows that the non-vanishing of the
amplitude given there is invariantly defined.

\medskip
\noindent
{\it Proof:}\;
By a partition of unity we can reduce the proof to the local case where 
$X_{j}\subset \mathbb{R}^{n_{j}}$, $T^{*}X_{j}$ is trivial for $j= 1, 2,$
$F=\int e^{\frac{i}{h}{\varphi}\left(x, z, \theta\right)} u\left(x, z, 
\theta\right)d\theta,$ 
where 
$\varphi$ is a non-degenerate phase function in a neighborhood of 
$\left(x_{0}, z_{0}, \theta_{0}\right)\in X_{1}\times X_{2}\times 
\mathbb{R}^{m}$ for some $m\in \mathbb{N}_{0}$ such that 
$\Lambda\cap U=\Lambda_{\varphi}$ for an open set $U$ with 
$\left(x_{0}, \xi_{0}; z_{0}, \eta_{0}\right)\in \Lambda\cap U,$ $u\in S_{n_{1}+n_{2}+m}^{\frac{n_{1}+n_{2}}{4}+\frac{m}{2}+r}(1)\cap C_{c}^{\infty}(\mathbb{R}^{n_{1}+n_{2}+m}),$ $u\sim\sum_{k=0}^{\infty}h^{k+\frac{n_{1}+n_{2}}{4}+\frac{m}{2}+r}u_{k},$  $u_{k}\in S_{n_{1}+n_{2}+m}^{0}(1),$ and
$A=\int_{\mathbb{R}^{n_{1}}} e^{\frac{i}{h}\left\langle x-y, \xi\right\rangle}a\left(x, 
\xi\right)d\xi$ with $a\in S_{2n_{1}}^{0}(1)\cap C_{c}^{\infty}(\mathbb{R}^{2n_{1}}),$ $a\sim\sum_{k=0}^{\infty}h^{k}a_k,$ $a_k\in S_{2 n_{1}}^{0}(1).$
%Why is $u$ elliptic?

Let $\Phi\left(y, \xi; x, z, \theta\right)=\left\langle x-y, \xi 
\right\rangle +
\varphi\left(y, z, \theta\right).$ 
Then $\Phi$ has a critical point 
\[p_0\left(x, z,
\theta\right)=\left(y_0\left(x, z, \theta\right), \xi_0\left(x, z, 
\theta\right)\right)=\left(x, \varphi'_{x}\left(x, z,
\theta\right)\right).\]
The Hessian of $\Phi$ is
\begin{equation*}
\Phi''\left(y_0\left(x, z, \theta\right), \xi_0\left(x, z, \theta\right); 
x, z, \theta\right)=
\begin{bmatrix}
\varphi''_{xx}\left(x, z, \theta\right) & -I\\
-I & 0
\end{bmatrix},
\end{equation*}
and has determinant 1 and signature 0.

Let $\Psi\left(w, \eta; x, z, \theta\right)=\varphi\left(x, w, 
\theta\right)+\left\langle w-z,
\eta\right\rangle.$ 
Then $\Psi$ has a critical point 
\[q_0\left(x, z, 
\theta\right)=\left(w_0\left(x, z, \theta\right),
\eta_0\left(x, z, \theta\right)\right)=\left(z, -\varphi'_z\left(x, z, 
\theta\right)\right).\]
The Hessian of $\Psi$ is 
\begin{equation*}
\Psi''\left(w_0\left(x, z, \theta\right), \eta_0\left(x, z, \theta\right); 
x, z, 
\theta\right)=\begin{bmatrix}
\varphi''_{zz}\left(x, z, \theta\right) & I\\
I & 0
\end{bmatrix},
\end{equation*}
and has determinant 1 and signature 0.

We define
\begin{equation*}
\begin{aligned}
g_{p_{0}\left(x, z, \theta\right)}\left(p\right)= & \Phi\left(p; x, z,
\theta\right)-\Phi\left(p_o\left(x,
z, \theta\right); x, z, \theta\right)\\
 & -\frac{\left\langle \Phi''\left(p_0\left(x, z, \theta\right); x, z,
\theta\right)\left(p-p_0\left(x, z, \theta\right)\right), p-p_0\left(x, z, 
\theta\right)\right\rangle}{2}
\end{aligned}
\end{equation*}
and
\begin{equation*}
\begin{aligned}
f_{q_{0}\left(x, z, \theta\right)}\left(q\right)= & \Psi\left(q; x, z,
\theta\right)-\Psi\left(q_o\left(x,
z, \theta\right); x, z, \theta\right)\\
 & -\frac{\left\langle \Psi''\left(q_0\left(x, z, \theta\right); x, z,
\theta\right)\left(q-q_0\left(x, z, \theta\right)\right), q-q_0\left(x, z, 
\theta\right)\right\rangle}{2}.
\end{aligned}
\end{equation*}
For $j\in\mathbb{N}_{0}$ and $c\in C_{c}^{\infty}\left(\mathbb{R}^{2m}\right),$ set
\begin{equation*}
\begin{aligned}
\left(L^{l}_j \left(cu\right)\right)\left(x, z, \theta\right)=\sum_{\nu-\mu=j} \sum_{2\nu \geq
3\mu} & \frac {\left\langle
\left(\Phi''\left(p_{0}\left(x, z, \theta\right); x, z, 
\theta\right)\right)^{-1} D,
D \right\rangle ^{\nu}}{i^{j}2^{\nu}\mu!\nu!}\\
& \quad \left(g_{p_{0}\left(x, z, \theta\right)}^{\mu} c\left(x, \cdot\cdot\right) 
u\left(\cdot, z,
\theta\right)\right)\left(p_0\left(x, z, 
\theta\right)\right).
\end{aligned}
\end{equation*}
and 
\begin{equation*}
\begin{aligned}
\left(L^{r}_j \left(uc\right)\right)\left(x, z, \theta\right)=\sum_{\nu-\mu=j} \sum_{2\nu \geq
3\mu} & \frac{\left\langle
\left(\Psi''\left(q_{0}\left(x, z, \theta\right); x, z, 
\theta\right)\right)^{-1} D,
D \right\rangle ^{\nu}}
{i^{j}2^{\nu}\mu!\nu!}\\
 & \quad \left(f_{q_{0}\left(x, z, \theta\right)}^{\mu} c\left(\cdot, 
\cdot\cdot\right) u\left(x, \cdot,
\theta\right)\right)\left(q_0\left(x, z, 
\theta\right)\right).
\end{aligned}
\end{equation*}

Now, since $\pi_{2}|_{\Lambda}$ is an immersion, it follows from the 
Inverse 
Function Theorem, 
that there exists $\kappa\in C^{\infty}\left(T^{*}\mathbb{R}^{n_2}; 
T^{*}\mathbb{R}^{n_1}\times T^{*}\mathbb{R}^{n_2}\right),$ 
such that $\kappa\circ\pi_2|_{\Lambda}=\text{id}|_{\Lambda}.$ 
Let, now, $b_{0}\in C_{c}^{\infty}\left(\mathbb{R}^{2n_{2}}\right)$
be supported near $\rho_2$ and satisfy
\begin{equation*}
b_{0}=i^{n_{1}-n_{2}}\kappa^{*}\left(\pi_{1}|_{\Lambda}\right)^{*}a_{0}.
\end{equation*}
Then 
$i^{n_{2}}j_{\varphi}^{*}\left(\pi_{2}|_{\Lambda}\right)^{*}b_{0}
-i^{n_{1}}j_{\varphi}^{*}\left(\pi_{1}|_{\Lambda}\right)^{*}a_{0}$ 
vanishes on $C_{\varphi}$ and since 
$\varphi$ is a 
non-degenerate phase function, it follows that there exist $c_{j}^{0}\in 
C^{\infty}_{c}\left(\mathbb{R}^{n_{1}+n_{2}+m}\right),$ $j=1, \dots, m,$ 
such that 
\[i^{n_{1}}j_{\varphi}^{*}\left(\pi_{1}|_{\Lambda}\right)^{*}a_{0}-i^{n_{2}}j_{\varphi}^{*}\left(\pi_{2}|_{\Lambda}\right)^{*}b_{0}
=\sum_{j=1}^{m} \varphi'_{\theta_{j}} c_{j}^{0}.\]
For every $k>0$ we now choose $b_k\in C_{c}^{\infty}\left(\mathbb{R}^{2n_{2}}\right)$ in such a way that
\begin{equation*}
\left(i^{n_1}\sum_{\alpha+\beta+\gamma=k}L_{\alpha}^{l}\left(a_{\gamma} u_{\beta}\right)-i^{n_{2}}\sum_{\alpha+\beta+\gamma=k,}L_{\alpha}^{r}\left(u_{\beta}b_{\gamma}\right)-\sum_{l=1}^{m}D_{\theta_{l}}\left(c_{l}^{k-1}u_{0}\right)\right)_{\big|_{C_{\varphi}}} =0
\end{equation*}
and choose $c_{l}^{k}\in 
C_{c}^{\infty}\left(\mathbb{R}^{n_1+n_2+m}\right),$ $j=1, 
\dots, m,$ such that
\begin{equation*}
i^{n_1}\sum_{\alpha+\beta+\gamma=k}L_{\alpha}^{l}\left(a_{\gamma} u_{\beta}\right)-i^{n_{2}}\sum_{\alpha+\beta+\gamma=k,}L_{\alpha}^{r}\left(u_{\beta}b_{\gamma}\right)-\sum_{l=1}^{m}D_{\theta_{l}}\left(c_{l}^{k-1}u_{0}\right) = \sum_{l=1}^{m}\varphi'_{\theta_{l}}c_{l}^{k}u_{0}.
\end{equation*}
Lastly, let
\begin{equation}\label{basympt}
b\sim\sum_{j=0}^{\infty} h^j b_j.
\end{equation}

In the integrals
\begin{equation*}
L\left(x, z\right)=AF\left(x, z\right)=\frac{1}{\left(2\pi h\right)^{n_1}} 
\int\int e^{\frac{i}{h}\left\langle 
x-y, 
\xi\right\rangle} 
a\left(x, \xi\right) e^{\frac{i}{h}\varphi\left(y, z, 
\theta\right)}u\left(y, z, \theta\right) d\theta 
dy d\xi
\end{equation*}
and
\begin{equation*}
R\left(x, z\right)=FB\left(x, z\right)=\frac{1}{\left(2\pi 
h\right)^{n_2}}\int\int e^{\frac{i}{h}\varphi\left(x, 
w,
\theta\right)}u\left(x, w, \theta\right) e^{\frac{i}{h}\left\langle w-z, 
\eta\right\rangle} b\left(w,
\eta\right)d\theta dw d\eta.
\end{equation*}
we now apply the method of stationary phase, Theorem 7.7.5 in \cite{H}, in the 
$\left(y, \xi\right)$ and the $\left(w, \eta\right)$
variables, respectively, and obtain 
\begin{equation*}
L\left(x, z\right) \sim i^{n_1}\sum_{t=0}^{\infty} h^{t-\frac{n_{1}+n_{2}}{4}-\frac{m}{2}-r}\sum_{j=0}^{t}\sum_{v=0}^{t-j} \int e^{\frac{i}{h}\varphi\left(x, z, \theta\right)} \left(L^{l}_{j} 
\left(a_{t-j-v}u_{v}\right)\right)\left(x, z, \theta\right) 
d\theta
\end{equation*}
\begin{equation*}
R\left(x, z\right) \sim i^{n_2}
\sum_{t=0}^{\infty} h^{t-\frac{n_{1}+n_{2}}{4}-\frac{m}{2}-r} \sum_{j=0}^{t}\sum_{v=0}^{t-j}\int e^{\frac{i}{h}\varphi\left(x, z, \theta\right)} 
\left(L^{r}_{j} \left(u_{v}b_{t-j-v}\right)\right)\left(x, z, \theta\right) d\theta.
\end{equation*}
By the choice of $b_0$ we have
\begin{equation*}
\begin{aligned}
 & \int e^{\frac{i}{h}\varphi\left(x, z, \theta\right)}u_{0}\left(x, z, 
\theta\right)\left[i^{n_{1}}a_{0}\left(x, \varphi'_{x}\left(x, z, \theta\right)\right)-i^{n_{2}}b_0\left(z, -\varphi'_{z}\left(x, z, \theta\right)\right)\right]d\theta\\
& \quad =\int e^{\frac{i}{h}\varphi\left(x, z, 
\theta\right)}\sum_{j=1}^{m}\varphi'_{\theta_{j}}\left(x, z, 
\theta\right)c_{j}^{0}\left(x, 
z, \theta\right)u_{0}\left(x, z, \theta\right)d\theta\\
& \quad =-h\int e^{\frac{i}{h}\varphi\left(x, z, 
\theta\right)}\sum_{j=1}^{m}D_{\theta_j}\left(c_{j}^{0}\left(x, 
z, \cdot\right)u_{0}\left(x, z, 
\cdot\right)\right)\left(\theta\right)d\theta,
\end{aligned}
\end{equation*}
using integration by parts.
By the choice of the symbol $b_1$ we then have that
\begin{equation*}
AF-FB=\mathcal{O}_{C_{c}^{\infty}(\mathbb{R}^{n_{1}+n_{2}})}\left(h^{2-\frac{n_{1}+n_{2}}{4}-\frac{m}{2}-r}\right).
\end{equation*}
Iterating this argument, we obtain from the choice of the symbols 
$b_k'$s that
\begin{equation*}
AF\equiv FB \text{ near } (\rho_1, \rho_2).
\end{equation*}

\noindent
{\bf Acknowledgements.}  I would like to thank Victor Ivrii and Maciej Zworski for helpful discussions during the preparation of this article.
I would also like to thank Vesselin Petkov for introducing me to the work of his student Laurent Michel, which has helped me complete my work on this article.


\begin{thebibliography}{99}
\bibitem{ChG}
C. Gerard,  Asymptotique des p\^{o}les de la matrice de 
scattering pour deux obstacles strictement convexes,  M\'emoires de la 
Soci\'et\'e Math\'ematique de France  116 (31) (1988).  
\bibitem{DS}
M. Dimassi and J. Sj\"{o}strand, Spectral Asymptotics in
the Semi-Classical Limit, Cambridge University Press, Cambridge, 1999.
\bibitem{GS}
A. Grigis and J. Sj\"{o}strand,  Microlocal Analysis for
Differential Operators, Cambridge University Press, Cambridge, 1994.
\bibitem{H}
L. H\"{o}rmander, The Analysis of Linear Partial Differential
Operators, Springer Verlag, Berlin, 1980.
\bibitem{M}
A. Martinez, An Introduction to Semiclassical and Microlocal
Analysis, Springer-Verlag, New York, 2002.
\bibitem{Michel}
L. Michel, Semi-classical Behavior of the Scattering Amplitude
for Trapping
Perturbations at Fixed Energy, Ph. D. Dissertation, Universit\'{e} de Bordeaux, 2001.
\bibitem{SZQ}
J. Sj\"{o}strand and M. Zworski, Quantum Monodromy and
Semi-classical Trace Formulae, Journal de Math\'{e}matiques Pures
at Appliqu\'{e}es  81 (1)(2002), 1--33.
\end{thebibliography}
\end{document}